\documentclass{amsart}

\usepackage{amsmath,amsfonts,amscd,amssymb,epsf}

\newtheorem{theorem}{Theorem}[section]
\newtheorem{lemma}[theorem]{Lemma}
\newtheorem{proposition}[theorem]{Proposition}

\newtheorem{corollary}[theorem]{Corollary}

\theoremstyle{definition}

\newtheorem{example}[theorem]{Example}

\theoremstyle{remark}
\newtheorem{remark}[theorem]{Remark}

\numberwithin{equation}{section}

\newcommand{\Hom}{\ensuremath{H\!om}}
\newcommand{\R}{\ensuremath{\mathbb R}}
\newcommand{\C}{\ensuremath{\mathbb C}}
\newcommand{\Z}{\ensuremath{\mathbb Z}}

\newcommand{\K}{\ensuremath{\mathbb K}}
\newcommand{\F}{\ensuremath{\mathbb F\,}}
\newcommand{\Q}{\ensuremath{\mathbb Q}}

\newcommand{\lqs}{{\lq\lq}}

\newcommand{\MCG}{\ensuremath{MCG\,\!}}
\newcommand{\PSL}{\ensuremath{{\mathbb P}SL\,\!}}
\newcommand{\SL}{\ensuremath{SL\,\!}}

\newcommand{\GL}{\ensuremath{GL\,\!}}
\newcommand{\SU}{\ensuremath{SU\,\!}}

\newcommand{\Out}{\ensuremath{Out\,\!}}
\newcommand{\Aut}{\ensuremath{Aut\,\!}}

\begin{document}

\author{Richard Brown}
\address{Department of Mathematics and Statistics,
The American University, 4400 Massachusetts Ave, NW Washington, DC 20016-8050 USA} \email{brown@american.edu}
\date{\today}


\title[Automorphisms of Fricke characters]
{Automorphisms of the Fricke characters of free groups}


\begin{abstract}
In this note, we embed the set of all Fricke characters of a free
group $F_n$ -- the set of all characters of representations of
$F_n$ into $\SL(\C^2)$ -- as an irreducible affine variety
$V_{F_n}\in\C^{2^n-1}$.  Using the Horowitz generating set as the
indeterminates, we show that the ideal $I_n$ of all polynomials in
these indeterminates which vanish on $V_{F_n}$ is generated by the
Magnus relation for arbitrary octets of elements in $\SL(\C^2)$.
Using this relation, we produce a basis for $I_n$, and show that
it is prime. We then show that the natural action of automorphisms
of $F_n$ on $V_{F_n}$ extends to polynomial automorphisms on all
of $\C^{2^n-1}$. Furthermore, for $\nu$ a complex volume (a
nonvanishing holomorphic $2^n-1$ form) on $\C^{2^n-1}$,
automorphisms of $F_n$ preserve $|\nu |$. This construction
provides an algebraic model for the analysis of the dynamics of
the measure preserving action of $\Out(F_n)$ on $V_{F_n}$.
\end{abstract}


\maketitle
\vskip .15truein


\section{Introduction}\label{sec:Intro}

In this paper, we construct a working algebraic model of the set of all special linear characters of the free
group on $n$-letters $F_n$.  This set arises as an affine variety $V_{F_n}\subset\C^{2^n-1}$.  Using the Horowitz
generating set (the characters of the $2^n-1$ basic words in $F_n$, described below) as the indeterminates for the
polynomial ring with integer coefficients, we extend the work of Magnus to construct finite a set of polynomials
that generate the ideal that defines $V_{F_n}$. $V_{F_n}$ can also be interpreted as the algebro-geometric
quotient of the set of all representation $F_n$ into $SL(\C^2)$, by the conjugacy map:  two representations are
identified if they are conjugate or if their conjugacy classes are inseparable as points of the usual quotient
(the term \lqs algebro-geometric" in this sense refers to the fact that the ring of functions on $V_{F_n}$
identifies with the ring of conjugation-invariant regular functions on the space of all representations).
$V_{F_n}$ is then called the $\SL(\C^2)$-character variety of $F_n$. The geometry of $V_{F_n}$ has been
extensively studied, and it is well known that automorphisms of $F_n$ act on $V_{F_n}$, preserving much of this
geometry. Here we show that under this Horowitz embedding, $V_{F_n}\subset\C^{2^n-1}$, the $\Out(F_n)$ action on
$V_{F_n}$ extends naturally to all of $\C^{2^n-1}$ as polynomial automorphisms which leave invariant the modulus
of the (complex) volume form formed by the exterior product of the differentials of the coordinates of
$\C^{2^n-1}$ given by the Horowitz generators.  We will refer to this form as the standard volume form.

\begin{theorem}\label{thm:Main} $V_{F_n}$ may be embedded in $\C^{2^n-1}$ as an irreducible affine variety
via the Horowitz generating set.  $\sigma\in\Out(F_n)$ induces $\widehat{\sigma}\in\Aut(V_{F_n})$ which extends to
a polynomial automorphism of $\C^{2^n-1}$ which preserves the modulus of the standard complex volume.\end{theorem}

A special case arises if we restrict to the real points of $V_{F_n}$, denoted $V_{F_n,\R}$.  Here
$V_{F_n,\R}\subset\R^{2^n-1}$ is a real affine variety. Then the real standard $2^n-1$ form is a true volume form.
If $F_n$ is a surface group -- the fundamental group of a compact orientable surface $S$ (necessarily with
boundary since $F_n$ is free) -- then $\SL(\C^2)$-characters of $F_n$ in two real forms \cite{MS}:
$\SL(\R^2)$-characters, and $\SU(2)$-characters ($2$-dimensional special unitary representations necessarily have
real characters). Depending on the genus of $S$ and the number of boundary components, $V_{F_n,\R}$ possesses the
additional structure of a Poisson space, with symplectic leaves corresponding to the characters of representations
which agree on the boundary components (see Huebschmann~\cite{Hueb}). The mapping class group of the surface
$\MCG(S)$ is in general a proper subset of $\Out(F_n)$ (it is the set of isotopy classes of orientation preserving
homeomorphisms of $S$ which pointwise fix the boundary of $S$. This condition imposes restrictions on the type of
automorphisms of $\pi_1(S)=F_n$ that correspond to mapping classes). In this case,
$\sigma\in\MCG(S)\subset\Out(F_n)$ preserves this Poisson structure, and acts symplectically on the leaves.
Mapping classes are necessarily volume preserving on these leaves, and hence on all of $V_{F_n,\R}$.

\begin{theorem}\label{thm:SubMain} Let $F_n=\pi_1(S)$ for $S$ a compact surface with boundary.  Then $\MCG(S)$ acts as
volume preserving automorphisms of $V_{F_n,\R}\subset\R^{2^n-1}$ which extend to volume preserving polynomial
automorphisms of $\R^{2^n-1}$.\end{theorem}

One can build a \lqs Poisson" volume form on $V_{F_n,\R}$ via the exterior product of the symplectic volume on the
leaves with a pull back of a volume form on the leaf space.  It would be interesting to know how the standard
volume form on $\R^{2^n-1}$ relates to the Poisson volume form on $V_{F_n,\R}$ under this Horowitz embedding.

It has also long been of interest to use $V_{F_n}$ or its associated ideal $I_n$ to study the properties of the
automorphism group $\Out(F_n)$ (see, for instance \cite{Hor},\cite{Wh},\cite{Mag}, \cite{Mc}, or \cite{GAMA}).  A
standard issue complicating this endeavor, however, has been to write a basis for $I_n$ when $n>3$.  Whittemore
produced a set of $6$ polynomials in $I_4$, yet could not show that this set generated all of $I_4$.  Magnus
(compare also \cite{GAMA}) recognized the utility of a \lqs general" identity (which we call the Magnus Relation
in Section~\ref{sec:Magnus}) for arbitrary octets of elements in $\SL(\C^2)$.  It turns out that all of the
members of $I_n$ can be derived via appropriate substitutions of matrices into the Magnus Relation. However,
Magnus' interest in finding an embedding of $V_{F_n}$ via a set of rational functions precluded him from actually
forming a proper implicit representation of $V_{F_n}$ as an affine set. It should be noted here that
Gonz\'alez-Acu\~na and Montesinos-Amilibia~\cite{GAMA} found a finite generating set for $I_n$ using a particular
subset of the Horowitz generators (see Remark~\ref{rem:GAMA}).  It turns out there that this subset does not lead
to volume preserving polynomial automorphisms of the ambient affine space (see Example~\ref{ex:GAMA}).  For the
sake of completeness in announcing the results of this paper, we also offer the following:

\begin{theorem}\label{thm:SubMain2} In the polynomial ring $\Z[x_1,\ldots,x_{2^n-1}]$ whose indeterminates are the
Horowitz generators of $F_n$, the ideal of polynomials that vanish identically for all characters of special
linear representations of $F_n$ is prime and finitely generated by $2^n-1-(3n-3)$ appropriate substitutions into
the Magnus Relation.\end{theorem}

This note is organized as follows:  In Section~\ref{sec:Horowitz}, we define the Horowitz generating set, and use
it to show that $V_{F_n}$ is an algebraic subset of $\C^{2^n-1}$. Section~\ref{sec:Magnus} deals with the Magnus
Relation, a general identity that we will use to create membership in the ideal which defines $V_{F_n}$.  In
Section~\ref{sec:CharVar}, we develop criteria for the construction of the generating polynomials for $I_n$
specified in Theorem~\ref{thm:SubMain2}, and show that $I_n$ is prime, thus establishing
Theorem~\ref{thm:SubMain2}. Using a standard generating set for $\Out(F_n)$ by Nielsen transformations, we will
show in Section~\ref{sec:OutFn} that these generators all act on $\C^{2^n-1}$ as polynomial automorphisms that
preserve $|\nu |$, where $\nu$ is a complex volume form on $\C^{2^n-1}$, thus proving Theorem~\ref{thm:Main}.
Theorem~\ref{thm:SubMain} will then follow once the first result is interpreted correctly. In
Section~\ref{sec:Proofs}, we recount the three main Theorems to consolidate their proofs.  The last section,
Section~\ref{sec:Examples}, presents some examples for $n=2$, $3$, and $4$.

\section{Horowitz Generators}\label{sec:Horowitz}

Fix an ordered generating set for $F_n$.  Call a word $X\in F_n$ {\em basic} if each letter of $X$ is a generator
of exponent one, and each letter of $X$ is greater than the one after it.  It is easy to see that there are
exactly $2^n-1$ basic words in $F_n$.  Denote this basic set by $\{X_i\}_{i=1}^{2^n-1}$, and extend the ordering
of the generating set of $F_n$ to $\{X_i\}_{i=1}^{2^n-1}$.  Indeed, let $wl:F_n\to\Z_+$ be the word-length
function (typically defined on cyclically reduced words in $F_n$.  The basic set is, however, cyclically reduced).
Extend the ordering to the basic set as follows:
\begin{itemize}
\item If $wl(X_i)>wl(X_j)$, then $X_j>X_i$ (shorter words are greater than longer words). \item If
$wl(X_i)=wl(X_j)$, then the ordering is determined by the ordering of the generating set at the first position
where the letters of $X_i$ and $X_j$ disagree.
\end{itemize}

\begin{example}\label{ex:basic} Let $F_3=\langle A,B,C\rangle$ ordered such that $A>B>C$.
Then the basic set for $F_3$ can be ordered from highest to lowest \[\{A,B,C,AB,AC,BC,ABC\}.\]\end{example}

For $\K$ a field, a special linear character of $X$ (the Fricke character of $X$) is defined via a representation
of $F_n$ into the Lie group $\SL(\K^2)$.  The character of a representation is the trace of the associated image
of $F_n$. By fixing the word $X\in F_n$, and varying the assigned matrix it is associated to, one obtains the
character of the word $X$. Thus the character of $X$ is a $\K$-valued, conjugation invariant map on the set of all
$\SL(\K^2)$-representations of $F_n$.  We will denote this map \[ tr_X:\Hom(F_n,\SL(\K^n))\to \K, \quad
tr_X(\phi)=tr(\phi(X)).\]  Adopt the notation that for a word $X\in F_n$, we will use an upper case letter, and
for its associated character, we will use the lower case equivalent.  Thus, $x$ is the character of the word
$X\in\F_n$ (that is, $x=tr_X$).

A result originally proposed by Fricke~\cite{Fr} was proved by Horowitz~\cite{Hor1}:

\begin{theorem}\label{thm:Hor} Given an arbitrary word $X\in F_n$, its character $x=tr_X$ can be expressed as a polynomial
with integer coefficients in the $2^n-1$ characters of the basic words in $F_n$.\end{theorem}

In this sense, the $2^n-1$ characters of the basic words, $\{ x_i\}_{i=1}^{2^n-1}$, generate the set of all
characters of $F_n$.  We will call this set the {\em Horowitz generating set}.  In Example~\ref{ex:basic}, the
Horowitz generating set of the characters of $F_3$ is
\[ \left\{x_i\right\}_{i=1}^7=\left\{a,b,c,ab,ac,bc,abc\right\}.\] Considering these generators as indeterminates,
the set of all characters of $F_n$ is then a subset of the ring of polynomials with integer coefficients in these
indeterminates $P_n$, where \[ P_n=\Z\left[ x_1,x_2,\ldots,x_{2^n-1}\right].\]

In general, however, there are polynomials with integer coefficients in these indeterminates which are identically
zero for any choice of representation.  These relations form an ideal $I_n$ within this polynomial ring, such that
the quotient ring $R_n=P_n/I_n$ is called the {\em ring of Fricke characters} (see~\cite{Mag}).  In
$P_{n,\C}=P_n\otimes\C$, the common zero locus of the elements in $I_{n,\C}$ is our definition of \[V_{F_n}=
R_{n,\C}=P_{n,\C}/I_{n,\C}.\]  Call $V_{F_n}$ the {\em Fricke variety} of $F_n$.

For $n>2$, $I_n$ is not trivial.  It was demonstrated by Fricke as well as many others, that for $n=3$, the Fricke
Relation generates a principal ideal $I_3$ in $P_3$. To see this, consider the following:

\begin{lemma}[Fricke]\label{lem:Fricke} Let $P$, $Q\in\Z[a,b,c,ab,ac,bc]$ denote the
polynomials \begin{eqnarray*} P&=& a\cdot bc+b\cdot ac+c\cdot ab -a\cdot b\cdot c\\
Q&=& a^2+b^2+c^2+ab^2+ac^2+bc^2+ab\cdot ac\cdot bc\\ &&-a\cdot b\cdot ab-a\cdot c\cdot ac-b\cdot c\cdot
bc-4.\end{eqnarray*}  Then $abc+acb = P$, $abc\cdot acb=Q$, and $abc$ and $acb$ are the roots of the quadratic
equation \[ z^2-Pz+Q=0.\]\end{lemma}

Based on this result, it is straightforward to see that the identity $abc^2-P\cdot abc+Q=0$ is a generator of the
principal ideal $I_3\subset P_3=\Z\left[ a,b,c,ab,ac,bc,abc\right]$: \[ I_3=\langle abc^2-P\cdot abc+Q\rangle.\]

For $n>3$, however, the number of generators of $I_n$ grows quickly, and a full generating set for $I_n$ has not
been explicitly given. Partial results have been constructed by Horowitz~\cite{Hor} and Whittemore~\cite{Wh} for
$n=4$. Using a general relation for special linear matrices, Magnus~\cite{Mag} found a way to embed $V_{F_n}$ via
a set of quadratic extensions (given by a modest set of relations) of the quotient field of the ring of
polynomials with integer coefficients in a subset of the Horowitz generators (those corresponding to the $3n-3$
words of either length 1, or length 2 with large first letter).  The remaining generators are then found via
rational functions of this minimal generating subset.  Below we will give a full implicit representation of
$V_{F_n}$ by constructing a basis of $I_n$.

\section{The Magnus Decomposition}\label{sec:Magnus}

In \cite{Mag} (compare also \cite{GAMA}), Magnus presents a \lqs general" identity, which holds for arbitrary
octets of matrices in $\SL(\C^2)$.  This relation is a generalization of the Fricke Relation above, and provides a
ready recipe for identifying new members of $I_n$. Let $M_\nu$, $N_\mu$, $\nu,\mu=1,2,3,4$ be any eight elements
of $\SL(\C^2)$. Denote by $tr(M_\nu N_\mu)$ the $4\times 4$-matrix, whose $ij$th element is the complex number
$tr(M_i N_j)$.

\begin{proposition}[Magnus Relation] $\det(tr M_\nu N_\mu)+\det(tr M_\nu N_\mu^{-1})=0.$\end{proposition}

Note that any representation of $F_n$ in $\SL(\C^2)$ is simply an assignment of the generators of $F_n$ to special
linear matrices.  Hence, the Magnus Relation is satisfied for any substitution of matrices assigned to a
representation;  any assignment of matrices to the Magnus relation results in a polynomial in $I_n$.  For this
discussion and the construction of $I_n$ in the next section, we will adopt the following notation:  For a group
element $X\in F_n$, we will also denote by $X$ its arbitrary assignment via a representation as a special linear
matrix.

\begin{example} For $F_3=\langle A,B,C\rangle$, let $M_1=N_1=A$, $M_2=N_2=B$, $M_3=N_3=AB$, and
$M_4=N_4=C$.  Then it is easy to show that we recover the Fricke Relation.  Indeed, the Magnus Relation becomes

\tiny
\[ \left| \begin{array}{rrrr}tr_{AA}&tr_{AB}&tr_{AAB}&tr_{AC}\\
tr_{BA}&tr_{BB}&tr_{BAB}&tr_{BC}\\ tr_{ABA}&tr_{ABB}&tr_{ABAB}&tr_{ABC}\\
tr_{CA}&tr_{CB}&tr_{CAB}&tr_{CC}\end{array}\right| + \left|
\begin{array}{rrrr}tr_{AA^{-1}}&tr_{AB^{-1}}&tr_{A(AB)^{-1}}&tr_{AC^{-1}}\\
tr_{BA^{-1}}&tr_{BB^{-1}}&tr_{B(AB)^{-1}}&tr_{BC^{-1}}\\ tr_{ABA^{-1}}&tr_{ABB^{-1}}&tr_{AB(AB)^{-1}}&tr_{ABC^{-1}}\\
tr_{CA^{-1}}&tr_{CB^{-1}}&tr_{C(AB)^{-1}}&tr_{CC^{-1}}\end{array}\right| = 0,\]

\normalsize which, upon utilizing some of the fundamental trace relations governing special linear characters (see
Fricke~\cite{Fr}) \[ tr_{X^{-1}}=tr_X, \] \[ tr_{XY}=tr_Xtr_Y-tr{XY^{-1}},\] reduces to

\tiny
\[ \left| \begin{array}{cccc}a^2-2&ab&a\cdot ab-b&ac\\
ab&b^2-2&b\cdot ab-a&bc\\ a\cdot ab-b&b\cdot ab-a&ab^2-2&abc\\
ac&bc&abc&c^2-2\end{array}\right| + \left|
\begin{array}{cccc}2&a\cdot b-ab&b&a\cdot c-ac\\
a\cdot b-ab&2&a&b\cdot c-bc\\ b&a&2&ab\cdot c-abc\\
a\cdot c-ac&b\cdot c-bc&ab\cdot c-abc&2\end{array}\right| = 0.\]

\normalsize

Writing out these determinants, this equation reduces to \[ (P-abc)\cdot abc-Q=0 \] where $P$ and $Q$ are
precisely as in Lemma~\ref{lem:Fricke}.
\end{example}

Indeed, for any basic word $X_i\in F_n=\langle A_1,\ldots,A_n\rangle$ of word-length $n>2$, $X_i$ may be divided
up into the concatenation of three other basic words in a number of ways in general:  \[ X_i=A_1A_2\cdots
A_n=(A_1\cdots A_j)\cdot(A_{j+1}\cdots A_k)\cdot(A_{k+1}\cdots A_n)=W_1\cdot W_2\cdot W_3.\] Then the assignments
$M_1=N_1=W_1$, $M_2=N_2=W_2$, $M_3=N_3=W_1\cdot W_2$, and $M_4=N_4=W_3$ into the Magnus Relation leads to an
identity in the corresponding Horowitz generators given by \[ z^2-Pz+Q=0, \] where \[ P, Q \in
Z[w_1,w_2,w_3,w_1w_2,w_1w_3,w_2w_3],\] and $z=w_1w_2w_3=a_1\cdots a_n$ is one of the roots.  It is worth noting
that the other root is $z=w_1w_3w_2=a_1\cdots a_j\cdot a_{k+1}\cdots a_n\cdot a_{j+1}\cdots a_k$.

Call this method of decomposing a basic word $X_i$ into a product of other basic words (of greater value in the
ordering) for the purpose of creating an identity using the Magnus Relation a {\em Magnus Decomposition} of $X_i$.
Some remarks are in order here:

\begin{remark} For a given $X_i$, where $wl(X_i)>3$, there are $\left(\begin{array}{c} wl(X_i)\\
3\end{array}\right)$ such Magnus Decompositions of $X_i$ of the above fashion alone (there are others also, as we
shall see).  It turns out that none are canonical, although there are criteria one can establish so that the
choice is always consistent.\end{remark}

\begin{remark} The above decomposition of a basic word into three other basic words of higher value in the
ordering is only one of many such types of decomposition.  Indeed, another which we shall employ, is to decompose
a basic word into two such pieces:  \[ X_i=A_1A_2\cdots A_n=(A_1\cdots A_\ell)\cdot (A_{\ell+1}\cdots
A_n)=W_1W_2.\]  The resulting assignments $M_1=N_1=A_1$, $M_2=N_2=A_2$, $M_3=N_3=W_1$, and $M_4=N_4=W_2$ results
in a polynomial which is either quadratic or linear, one of whose roots is $x_i$, and whose coefficients lie in
the polynomial ring generated by $a_1$, $a_2$, $w_1$, $w_2$, and the generators formed by the various pairings of
these.\end{remark}

\begin{remark} A Magnus Decomposition of a generator $x_i$ in the above fashion produces a polynomial in the
Horowitz generators, where $x_i$ is the generator os least value in the ordering the monomials inherit from the
ordering of the basic word set in $F_n$.  This fact will be central to the utility of the procedure we will employ
to construct the variety $V_{F_n}$.\end{remark}

\section{The character variety $V_{F_n}$}\label{sec:CharVar}

The following embedding of the ring of Fricke characters $R_n$ is given by Magnus~\cite{Mag}.  It is an implicit
representation of $R_n$ as the zero locus of a combination of polynomials and rational functions.  Recall that in
this context, the Fricke variety $V_{F_n}$ is simply the set of complex points of $R_n$:  $V_{F_n}=R_{n,\C}$ (See
Section~\ref{sec:Horowitz}):

\begin{theorem}  Let $\Omega_n$ be the quotient field of the ring of polynomials with integral coefficients in the
indeterminates $\{x_i\}_{i=1}^{3n-3}$ formed from the basic words of $F_n$ with word-length $1$ or $2$ with
initial generator $X_1$ or $X_2$.  Then $R_n$ is embeddable via a set of rational functions in an algebraic
extension of $\Omega_n$ which consists of at most $n-2$ simultaneous quadratic extensions.\end{theorem}

The theorem assumes that the dimension of the ambient space for the variety should be as small as possible.  We
prefer the following tack:

\begin{theorem}  Let $P_n$ be the ring of polynomials with integral coefficients in the
indeterminates $\{x_i\}_{i=1}^{2n-1}$ formed from the Horowitz generating set.  Then $V_{F_n}$ is the zero locus
in $P_n$ of the set of polynomials formed from all substitutions of characters of $F_n$ into the Magnus
Relation.\end{theorem}

Call the ideal consisting of all substitutions of characters of $F_n$ into the Magnus Relation $I_M$.  The
construction of the ideal $I_n$ whose variety is $V_{F_n}$ will consist of building a finite basis for $I_M$ and
showing that $I_n=I_M$.  We will follow the tack of Magnus in using the $3n-3$ initial indeterminates (those with
the largest values in the monomial ordering from Section~\ref{sec:Horowitz}) as the foundation set. Then for each
Horowitz generator $x_j$, $j>3n-3$, after this initial set, one can choose a particular set of characters for
substitution into the Magnus relation to generate a polynomial whose monomial set includes $x_j$ and no generators
lower in value than it.  This technique will result in a set of $2^n-1-(3n-3)$ distinct polynomials, whose ideal
is precisely $I_n$.

Indeed, for the following construction, consider the notation:  For the indeterminates $\{x_j\}_{j=3n-3}^{2^n-1}$,
denote $y_i=x_{3n-3+i}$.  Then $P_n=\Z [x_1,\ldots,x_{3n-3},y_1,\ldots,y_m]$, where $m=2^n-1-(3n-3)$.  The latter
set is distinguished from the former due to their constructive role of $I_n$.  For each Horowitz generator $y_i$,
denote the corresponding polynomial created $p_{y_i}$, and the corresponding basic word $Y_i\in F_n$.  According
to the previous section, there are many ways to decompose a basic word $Y_i$.  We will adopt the following
convention for the construction of $p_{y_i}$:

\vspace{8pt}
\begin{itemize}\partopsep 20pt \itemsep 10pt
\item[Case 1:] $wl(Y_i)=2$, where $Y_i=A_\mu A_\nu$ and $\nu>\mu>2$. \newline Choose $M_1=N_1=A_1$, $M_2=N_2=A_2$,
$M_3=N_3=A_\mu$ and $M_4=N_4=A_\nu$.  With these substitutions, the Magnus Relation produces a polynomial
$p_{y_i}$ which is either irreducible quadratic or linear, one of whose roots is $y_i$, and whose coefficients are
elements of the ring
\[ \Z[a_1,a_2,a_\mu,a_\nu,a_1a_2,a_1a_\mu,a_1a_\nu,a_2a_\mu,a_2a_\nu].\]  Note that the total degree of
$p_{y_i}$ will be 8 in this case.

\item[Case 2:] $wl(Y_i)>2$, where $Y_i=A_\mu A_\nu X_j$, and $wl(A_\mu)=wl(A_\nu)=1$. \newline Use the above
Magnus decomposition of $Y_i$ into three pieces and choose $M_1=N_1=A_\mu$, $M_2=N_2=A_\nu$, $M_3=N_3=A_\mu A_\nu$
and $M_4=N_4=X_j$.  Then $p_{y_i}$ will be an irreducible quadratic, one of whose roots is $y_i$ and whose
coefficients are elements in the ring \[ \Z[a_\mu,a_\nu,x_j,a_\mu a_\nu,a_\mu x_j,a_\nu x_j].\]  The total degree
of $p_{y_i}$ here is 4.
\end{itemize}

\subsection{The construction of $I_n$.}  Denote by $I_{n,0}=\langle 0\rangle \subset P_{n,0}=\Z[x_1,\ldots,x_{3n-3}]$.
Define $p_{y_1}=z^2-Pz+Q$ according to the above construction involving the Magnus Relation.  Note that
$p_{y_1}\in P_{n,0}$ is a monic with no roots in $P_{n,0}$.  However, one of the roots of $p_{y_1}$ is precisely
the new indeterminate $y_1$. Construct $P_{n,1}=P_{n,0}[y_1]$ as an integral extension.  Then \[
p_{y_1}(y_1)=y_1^2-Py_1+Q\] is a polynomial in $P_{y,1}$, whose value is identically 0 for all values of
$x_1,\ldots,x_{3n-1},y_1$ that correspond to characters of representations of $F_n$.  Thus $p_{y_1}(y_1)$
generates an ideal in $P_{n,1}$, which we will denote $I_{n,1}$.  Note that to avoid an excessive use of notation,
we will refer to the defining element in $I_{n,1}$ simply as $p_{y_1}$ (that is, $I_{n,1}=\langle
p_{y_1}\rangle$).

Inductively, suppose $1\le i<m$, and let $I_{n,i}=\langle p_{y_1},\ldots, p_{y_i}\rangle$ be an ideal in
$P_{n,i}=\Z[x_1,\ldots,x_{3n-3},y_1,\ldots,y_i]$.  Again, using the Magnus Relation construction above, define
$p_{y_{i+1}}=z^2-Pz+Q$, a monic in $P_{n,i}$, one of whose roots is $y_{i+1}$.  It is clear that
$p_{y_{i+1}}\not\in I_{n,i}$.  Extend $P_{n,i+1}=P_{n,i}[y_{i+1}]$, and construct a new ideal
$I_{n,i+1}=I_{n,i}+\langle p_{y_{i+1}}\rangle\in P_{n,i+1}$.

When $i=m$, $P_{n,m}=P_n$, and define compatibly $I_n=I_{n,m}=\langle p_{y_1},\ldots,p_{y_m}\rangle$.

\begin{lemma} The ideal $I_{n,i}=\langle p_{y_1},\ldots,p_{y_{i-1}},p_{y_i}\rangle$ is independent of the choice of
Magnus Decomposition of $Y_i\in F_n$ in the construction of $p_{y_i}$.\end{lemma}

\begin{proof} Construct $I_{n,i}=\langle p_{y_1},\ldots,p_{y_{i-1}},p_{y_i}\rangle$ and
$I_{n,i}^\prime=\langle p_{y_1},\ldots,p_{y_{i-1}},p_{y_i}^\prime\rangle$ based on two different choices of Magnus
Decomposition for the construction of the $y_i$ polynomial.  For any point
$x=(x_1,\ldots,x_{3n-3},y_1,\ldots,y_i)\in V(I_{n,i})$, we have
\[ p_{y_1}(x)=\cdots =p_{y_{i-1}}(x)=0.\]  By construction, both $p_{y_i}$ and $p_{y_i}^\prime$ are polynomials in
$P_{n,i}$ constructed precisely to have $y_i$ (the last coordinate of $x$) as a root.  Hence \[
p_{y_i}(x)=p_{y_i}^\prime(x)=0,\] and $x\in V(I_{n,i}^\prime)$.  By symmetry, the zero locus of both $I_{n,i}$ and
$I_{n,i}^\prime$ are the same.  Hence the ideals are the same.\end{proof}

For the following discussion, we will consider the ring $P_n$ to be the ring of polynomials with complex
coefficients $P_n\otimes\C$, and drop the subscript $\C$.

\begin{proposition}\label{thm:Prime} $I_n$ is prime in $P_n$.\end{proposition}

\begin{proof}  Since $\langle p_{y_1}\rangle=I_{n,1}$, and $p_{y_1}$ is irreducible, $I_{n,1}$ is prime in
$P_{n,1}$.  To proceed inductively, suppose $I_{n,i}$ is prime in $P_{n,i}$, where \[
P_{n,i}=\C[x_1,\ldots,x_{3n-3},y_1,\ldots,y_i].\]  Construct $p_{y_{i+1}}\in P_{n,i+1}=P_{n,i}[y_{i+1}]$ in the
above fashion, and note that $p_{y_{i+1}}$ is irreducible over $\C$.  It is sufficient to show that
$I_{n,i+1}=I_{n,i}+\langle p_{y_{i+1}}\rangle$ is prime in $P_{n,i+1}$.  To this end, note that since $I_{n,i}$ is
prime, it is radical, and $R_{n,i}=P_{n,i}/I_{n,i}$ is an integral domain (so there are no nonzero divisors).  In
$P_{n,i+1}$, view $I_{n,i}\subset I_{n,i+1}$, so that the mapping \[ P_{n,i+1}/I_{n,i+1}\longrightarrow
(P_{n,i+1}/I_{n,i})/(I_{n,i+1}/I_{n,i})\] is an isomorphism of rings.  Note that \[ P_{n,i+1}/I_{n,i}\cong
P_{n,i}[y_{i+1}]/I_{n,i}=R_{n,i}[y_{i+1}] \] since by elimination of $y_{i+1}$, $I_{n,i}\cap P_{n,i}=I_{n,i}$. And
$I_{n,i+1}/I_{n,i}\cong \langle p_{y_{i+1}}\rangle$, so that \[ P_{n,i+1}/I_{n,i+1}\cong R_{n,i}[y_{i+1}]/ \langle
p_{y_{i+1}}\rangle.\]

Let $f,g\in P_{n,i+1}$, such that $fg\in I_{n,i+1}$ (we will use brackets around polynomials to denote their
corresponding elements in $R_{n,i}[y_{i+1}]$).  Then $[fg]=0$ in $R_{n,i}[y_{i+1}]/ \langle p_{y_{i+1}}\rangle$.
Hence \[ [fg]=[k]p_{y_{i+1}}\in R_{n,i}[y_{i+1}] \] for some polynomial $[k]$.  Thus $p_{y_{i+1}}$ divides $[fg]$.
Since $p_{y_{i+1}}$ is irreducible over $\C$, $p_{y_{i+1}}$ must divide either $[f]$ or $[g]$.  Hence either
$[f]\in\langle p_{y_{i+1}}\rangle$ or $[g]\in\langle p_{y_{i+1}}\rangle$.  But then either $f\in I_{n,i+1}$ or
$g\in I_{n,i+1}$.  Hence $I_{n,i+1}$ is prime.\end{proof}

\begin{corollary}\label{thm:irr} $V(I_n)$ is irreducible.\end{corollary}

\begin{proof} As $I_{n,m}=I_n$ by definition, and $I_n$ is prime, this is obvious.\end{proof}

\begin{theorem}\label{thm:ideal} $V(I_{m,n})=V(I_n)=V_{F_n}$.\end{theorem}

To prove this theorem, we will make use of the following theorem:

\begin{theorem}[Magnus~\cite{Mag}]\label{thm:mag}  Let $A_1, A_2, A_3\in\SL(\C^2)$ such that their corresponding traces
\[ \begin{array}{rrr}tr_{A_1}=x_1,&tr_{A_2}=x_2,&tr_{A_3}=x_3,\\ tr_{A_1A_2}=x_{n+1},&tr_{A_1A_3}=x_{n+2},
&tr_{A_2A_3}=x_{2n}\end{array}\] satisfy the following 2 conditions:
\begin{itemize}\partopsep 20pt \itemsep 10pt
\item[1:] $tr_{[A_1,A_2]}=x_1^2+x_2^2+x_{n+1}^2-x_1x_2x_{n+1}-2\not= 2$, and

\item[2:] The discriminant of $p_{y_1}=z^2-Pz+Q$ defined above, where $y_1=tr_{A_1A_2A_3}$, doesn't vanish. That
is, $P^2-4Q\not=0$.
\end{itemize}

Then there exist matrices $A_4,\ldots,A_n\in\SL(\C^2)$ such that all of the rest of $\{x_i\}_{i=1}^{3n-3}$ assume
preassigned values.  At this point, then, the rest of the coordinates $\{y_i\}_1^m$ can be found by solving the
quadratic equations which generate $I_{n,m}$.\end{theorem}

\begin{proof}[Proof of Theorem~\ref{thm:ideal}] By construction, the character of any representation satisfies every manifestation of the Magnus
relation.  Denote the ideal generated by all manifestations of the Magnus Relation by $I_M$.  Then
\[ V_{F_n}\subset V(I_M)\subset V(I_{n,m}).\]  Since, a priori, these inclusions may be strict, it remains to show
that $V(I_{n,m})\subset V_{F_n}$.

To this end, let $x\in V(I_{n,m})$.  Then $x=(x_1,\ldots,x_{3n-3},y_1,\ldots,y_m)$, where $p_{y_i}(x)=0$, for
$i=1,\ldots, m$.  By Theorem~\ref{thm:mag}, if the first three coordinates of $x$ satisfy the two conditions of
the theorem, then there will be a representation of $F_n$ in $\SL(\C^2)$ whose character will assume the rest of
the $x_i$'s.  The corresponding $y_i$'s will then follow.  Hence, $x\in V_{F_n}$.  Note here that $V(I_{n,m})$ is
irreducible by Corollary~\ref{thm:irr}.

The two conditions above in the theorem are open conditions on both $\SL(\C^2)$ and $\F^{2^n-1}$, so that if they
hold anywhere on $V(I_{n,m})$, they will hold on all of at least the nonsingular part of $V(I_{n,m})$.
Since $V_{F_n}\subset V(I_{n,m})$, consider any representation where \[ A_i=\left(\begin{array}{cc} \sqrt{i}&\sqrt{i+1}\\
\frac{2}{\sqrt{i+1}}&\frac{3}{\sqrt{i}}\end{array}\right)\in\SL(\C^2)\quad i=1,2,3.\] It is easy to check that
$tr_{[A_1A_2]}\cong -0.5598\not= 2$ and that the discriminant of $p_{y_1}$, namely $P^2-4Q\cong 0.0239\not= 0$.

Hence ${\textsl Nonsing}(V(I_{n,m}))\subset V_{F_n}$.  By irreducibility, then, and the fact that $V_{F_n}$ is
closed, \[ V(I_{n,m})=\overline{{\textsl Nonsing}(V(I_{n,m}))}\subset V_{F_n}.\]\end{proof}

\begin{example} In general, $\dim_H(V_{F_n})=3n-3$, and this agrees with the interpretation of $V_{F_n}$
as the $\SL(\C^2)$-representation variety of $F_n$.  Indeed, the set of all representations of $F_n$ into
$\SL(\C^2)$ is simply the set of all $n$-tuples of matrices in $\SL(\C^2)$.  Hence \[ \Hom(F_n,\SL(\C^2)) =
\SL(\C^2)^n.\]  The complex dimension of $\Hom(F_n,\SL(\C^2))$ is $3n$.  Under the categorical quotient of
$\Hom(F_n,\SL(\C^2))$ by diagonal conjugation (mod out by the closures of orbits to eliminate the pathologies
associated to nonseparable orbits (see Goldman~\cite{Go} for an example)), the dimension of the quotient
$\Hom(F_n,\SL(\C^2))/\SL(\C^2)$ is $3n-3$.

For $n=4$, the representation used in the above proof of Theorem~\ref{thm:ideal}, along with the corresponding
assignment to $A_4$, leads to a point in $p\in\C^{2^4-1}=\C^{15}$ which is nonsingular in $V_{F_4}$. Indeed, it
can be verified through direct calculation that the $6\times 15$ matrix of partial derivatives of the generators
of $I_{6,4}$ (that is, the Jacobian
\[J_p(p_{y_1},\ldots,p_{y_6})=\left(\frac{\partial p_{y_i}}{\partial x_j}\right)\] has rank $6$.  Thus the
dimension of $V_{F_4}$ is $15-6=9$.  Direct calculation also reveals that this is true even if $p_{y_6}$, which
corresponds to a particular Magnus decomposition of the basic word $A_1A_2A_3A_4$ is replaced by any other of the
4 choices.  Note that the choice of Magnus decomposition used to create $p_{y_1}$ was $A_1\cdot A_2\cdot A_3A_4$.
The other choices are $A_1\cdot A_2A_3\cdot A_4$, $A_1A_2\cdot A_3\cdot A_4$, and $A_1A_4\cdot A_2\cdot
A_3$.\end{example}

\begin{remark} $I_n$ defined herein as $I(V_{F_n})$ is the ideal studied by Horowitz~\cite{Hor} and
Whittemore~\cite{Wh}.  In fact, Whittemore's partial basis of $I_4$ includes the polynomial $p_5$.  A close
inspection reveals that this polynomial corresponds to the Magnus decomposition of the basic word $A_1A_2A_3A_4\in
F_4$ given by $A_1A_2\cdot A_3\cdot A_4$.  She then shows that the other generators are obtained as the images of
this (and hence of the sums of images of those already found) under suitable automorphisms of $F_n$. It is obvious
that automorphisms of $F_n$ preserve $I_n$ (they simply take characters to characters), and in many cases either
simply permute the Magnus generators of $I_n$, or take a particular generator corresponding to one choice of a
Magnus decomposition of a basic word, into that of another.\end{remark}

\begin{remark}\label{rem:GAMA} Gonz\'alez-Acu\~na and Montesinos-Amilibia~\cite{GAMA} also find a finite basis for
the ideal of polynomials that vanish for all $\SL(\C^2)$-characters of $F_n$.  In their construction, they choose
the set of all Horowitz generators whose basic words are of length three or less.  This is a set of $p$
indeterminates, where $p=\frac{n(n^2+5)}{6}$.  For $n=4$, this set is the $14$ elements \[
\{a,b,c,d,ab,ac,ad,bc,bd,abc,abd,acd,bcd\}.\]  The Horowitz generators of higher word length are then found to be
polynomials in the Horowitz generators of word length three or less. This, like the construction in
Magnus~\cite{Mag}, provides an embedding of $V_{F_n}$ in a smaller ambient affine space. However, the ideal $I_n$
is not in $P_n$, but rather $P_n\otimes\Q$, and the natural extension of $\Out(F_n)$ to maps of $\C^p$ is not
necessarily by polynomial automorphisms (in particular, their determinants are not nonvanishing). See
Example~\ref{ex:GAMA}.\end{remark}

\section{Polynomial automorphisms of $\C^{2^n-1}$}\label{sec:OutFn}

Recall that $P_{n,\C}=P_n\otimes \C=\C [x_1,\ldots ,x_{2^n-1}]$ be the polynomial ring on $2^n-1$ indeterminates
with complex coefficients. By considering the indeterminates of $P_n$ as the coordinates of $\C^{2^n-1}$, it is
easy to see that $\sigma\in\Aut(P_{n,\C})$ induces a polynomial automorphism $\widehat{\sigma}\in
PolyAut(\C^{2^n-1})$ (It induces a polynomial map.  It is an automorphism since the inverse automorphism of
$P_{n,\C}$ also exists and induces an inverse to $\widehat{\sigma}$ which is also polynomial). Let
$\sigma\in\Aut(P_n)$.  By Whittemore~\cite{Wh}, $\Out(F_n)\subset\Aut_{I_n}(P_n)\subset\Aut(P_n)$.  Extend this to
a group of automorphisms $\Aut(P_{n,\C})$.  In this section, we will show that $\sigma\in\Out(F_n)$ induces a
$\widehat{\sigma}$ such that $\det( Jac(\widehat{\sigma}))\equiv \pm 1$:

For $F_n=\langle A_1,\ldots,A_n\rangle$, Nielsen~\cite{Niel} presents a generating set for $\Out(F_n)$ given by
the four elements below, named respectively, twist, two-element permutation, cyclic permutation, and inversion:
\[ \Phi_1:A_1\mapsto A_1A_2 \quad \Phi_2:\begin{array}{rcl} A_1&\mapsto&A_2 \\
A_2&\mapsto&A_1\end{array}\quad \Phi_3:\begin{array} {rcl} A_i&\mapsto&A_{i+1}\\ A_n&\mapsto&A_1\end{array} \quad
\Phi_4:A_1\mapsto A_1^{-1}.\]

\begin{proposition}\label{thm:Jac1} For any $\sigma\in\Out(F_n)$, the induced map $\widehat\sigma\in PolyAut(\C^{2^n-1})$
satisfies $\det\left( Jac(\widehat\sigma)\right)\in\{+1,-1\}$.\end{proposition}

\begin{proof} The above generators of $\Out(F_n)$ are all of finite order except for $\Phi_1$.  Hence, any
$\sigma\in\Out(F_n)$ in the subgroup generated by $\Phi_2$, $\Phi_3$, and $\Phi_4$ will necessarily induce a
finite order element $\widehat\sigma\in PolyAut(\C^{2^n-1})$.  Thus some iterate of $\widehat\sigma$ must be the
identity, and hence the determinate of the Jacobian of this iterate of $\widehat\sigma$ must be everywhere $1$.
Thus the determinate of the Jacobian of $\widehat\sigma$ must be a constant root of 1.  As $\Phi_2$ and $\Phi_4$
are involutions, their induced maps must be also, hence the proposition is satisfied for these two generators.
$\Phi_3$ is of order $n$.  However, the induced map is simply a permutation of the Horowitz generators, and hence
the Jacobian will be constant and an element of $\GL(\Z^{2^n-1})$.  Hence, its determinate will either be a
constant $1$ or $-1$.  The theorem will be proved once we establish the result for $\Phi_1$.  This is the content
of Lemma~\ref{thm:Phi1} below. \end{proof}

\begin{remark} The finite order generators of $\Out(F_n)$, namely $\Phi_2$, $\Phi_3$, and $\Phi_4$, form an
interesting subgroup.  The subgroup $\langle \Phi_2, \Phi_3\rangle$ is simply the full symmetric group on
$n$-letters (full permutation group of the generators of $F_n$).  With $\Phi_4$, the group $\langle \Phi_2,
\Phi_3,\Phi_4\rangle$ is just the group of signed permutations, or the hyperoctahedral group of order $n$ (see
McCool~\cite{Mc}).  This group has $2^nn!$ elements.
\end{remark}

\begin{remark} That the $\Out(F_n)$ action on $\C^{2^n-1}$ is by polynomial maps is an easy consequence of Horowitz's
Theorem~\ref{thm:Hor}.  Showing the polynomial maps generated by the $\{\Phi_i\}_{i=1}^4$ are invertible is also
straightforward, since the inverses are easily constructible and are polynomial.  In fact, it is readily apparent
that $\widehat\Phi_1$ and $\widehat\Phi_4$ are quadratic (the maximal degree of the component polynomials is $2$),
and $\widehat\Phi_3$ is linear (it is a simple permutation of the coordinates).  It is known, see
Wang~\cite{Wang}, that the Jacobian Conjecture holds in degrees $1$ and $2$:  A quadratic or linear polynomial map
on $\C^n$ with an everywhere nonvanishing Jacobian is necessarily an automorphism. Hence, except for
$\widehat\Phi_2$, showing the Jacobian of each $\Phi_i$ doesn't vanish is sufficient to establish that they are
automorphisms. $\widehat\Phi_2$, however, is in general a cubic map, and the Jacobian Conjecture is still a
conjecture in higher degree. Since $\widehat\Phi_2$ is an involution, though, its square is the identity. Hence
the inverse of $\widehat\Phi_2$ is polynomial, and $\widehat\Phi_2$ is an automorphism.\end{remark}

\begin{lemma}\label{thm:Phi1} $\det\left( Jac(\widehat\Phi_1)\right) = 1$.\end{lemma}

To prove this, return to the notation where $F_n=\langle A, B, C,\ldots\rangle$ is ordered so that $A>B>C>\ldots$,
and $\{x_i\}$ are the $2^n-1$ Horowitz generators of $P_n$ are ordered as stipulated in
Section~\ref{sec:Horowitz}. The Horowitz generators of $F_3$, for example, are the 7 functions \[ x_1=a,\; b,\;
c,\;ab,\; ac\; ,bc\; ,abc = x_7.\]

\begin{lemma} For $n>2$, there exist an even number of Horowitz generators of the form $abv$, where $v$ is comprised of
other generators.\label{lem:hor}\end{lemma}

\begin{proof} $abv$ corresponds to the basic word $ABV\in F_n$.  Thus $V$ is basic, and $A,B\not\in V$.
For $n>2$, there are $2^{n-2}-1$ such basic words $V$, which is odd.  With $V=e$ as another choice, the number is
now even.\end{proof}

Let $r_i$ be the $i$th row of $Jac(\widehat\Phi_1)$.  Then $r_i=\left( \frac{\partial}{\partial x_j}
\widehat\Phi_1(x_i)\right)$. Denote by $\lfloor\widehat\Phi_1(x_i)\rfloor$ the index of the last non-zero element
in $r_i$.

\begin{lemma}\label{lem:ri} $\forall x_i\not= av,$ where $b\not\in v$,
$\lfloor\widehat\Phi_1(x_i)\rfloor= i$.\end{lemma}

\begin{proof} For any word $W\in F_n$ such that $A\not\in W$, $\Phi_1(W)=W$.  Hence $\widehat\Phi_1(x_i)=x_i$
for any Horowitz generator where $a\notin x_i$.  Thus in this case, $\lfloor\Phi_1(x_i)\rfloor=i$.  If $x_i$ is of
the form $abv$, for some word $v$, which is either basic or $e$, then \begin{equation*} \widehat\Phi_1(abv) = abbv
= abv\cdot b-av.\end{equation*} Thus in this case also, $\lfloor\widehat\Phi_1(x_i)\rfloor =i$. This exhausts the
supply of Horowitz generators stipulated in the lemma.\end{proof}

In fact, the lemma does not hold only in the cases where \begin{equation*}
x_i=av\stackrel{\widehat\Phi_1}{\longmapsto} abv.\end{equation*}

\begin{proof}[proof of Lemma~\ref{thm:Phi1} above] By Lemma~\ref{lem:ri}, for all $x_i$ not of the form $av$,
$\lfloor\widehat\Phi_1(x_i)\rfloor=i$.  Thus, $Jac(\widehat\Phi_1)$ is almost lower triangular.  For each $v$ a
Horowitz generator, such that $a,b\not\in v$, there is a generator pair of the form $\{av,abv\}$.  Under the
action by $\widehat\Phi_1$, \begin{equation*} \begin{array}{rcl} \widehat\Phi_1(av)&=&abv\\ \widehat\Phi_1(abv)&=&
abv\cdot b-av.\end{array} \end{equation*}

Let $e_v$ correspond to the elementary column operation on matrices which is the switch of the two columns
corresponding to $av$ and $abv$.  Note after performing this switch, \begin{equation*} \begin{array} {rclcl}
av&\stackrel{\widehat\Phi_1}{\longmapsto}&abv&\stackrel{e_v}{\longmapsto}&av\\ abv&\longmapsto&abv\cdot
b-av&\longmapsto&av\cdot b-abv.\end{array}\end{equation*}

This operation affects no other rows of $Jac(\widehat\Phi_1)$ other than the rows corresponding to $av$ and $abv$.
For each switch, \begin{equation*} \det(e_v Jac(\widehat\Phi_1)) = -\det(Jac(\widehat\Phi_1)).\end{equation*}
However, by Lemma~\ref{lem:hor}, there are an even number of pairs. Hence after all of the switches, the
determinants will be the same.

After all of these switches, the transformed Jacobian will indeed be lower triangular.  Hence, its determinant can
be found by simply multiplying all of the elements in the main diagonal.  For all $x_i\not= abv$, the element on
the main diagonal is 1.  For each $x_i=abv$, the element on the main diagonal is $-1$.  Since there are en even
number of them, it follows that $\det(Jac(\widehat\Phi_1))=1$.\end{proof}

\section{Proofs of theorems}\label{sec:Proofs}

The three theorems mentioned in Section~\ref{sec:Intro} are now consequences of the constructions of $I_n$ and
$V_{F_n}$.  In this section, we consolidate the content of the previous sections and address the proofs of the
three theorems on Section~\ref{sec:Intro} directly. We start with Theorem~\ref{thm:SubMain2}:

\setcounter{section}{1} \setcounter{theorem}{2}
\begin{theorem} In the polynomial ring $\Z[x_1,\ldots,x_{2^n-1}]$ whose indeterminates are the
Horowitz generators of $F_n$, the ideal of polynomials that vanish identically for all characters of special
linear representations of $F_n$ is prime and finitely generated by $2^n-1-(3n-3)$ appropriate substitutions into
the Magnus Relation.\end{theorem}

\begin{proof} In Section~\ref{sec:CharVar}, we establish that $I_n$ is precisely the ideal of all polynomials that
vanish identically for all characters of $\SL(\C^2)$-representations $V_{F_n}$.  By Proposition~\ref{thm:Prime},
$I_n$ is prime, and Theorem~\ref{thm:ideal} shows that $V(I_n)$ is precisely $V_{F_n}$, when $I_n$ is generated by
the $2^n-1-(3n-3)$ polynomials constructed via substitutions into the Magnus Relation given by appropriate Magnus
decompositions of the basic words in $F_n$.\end{proof}

We are now in position to prove Theorem~\ref{thm:Main}:

\setcounter{theorem}{0}
\begin{theorem}  $V_{F_n}$ may be embedded in $\C^{2^n-1}$ as an irreducible affine variety
via the Horowitz generating set.  $\sigma\in\Out(F_n)$ induces $\widehat{\sigma}\in\Aut(V_{F_n})$ which extends to
a polynomial automorphism of $\C^{2^n-1}$ which preserves the modulus of the standard complex volume
form.\end{theorem}

\begin{proof} The fact that $V_{F_n}$ is an affine variety is given by Theorem~\ref{thm:ideal}.  It is irreducible
since $I_n$ is prime.  Any element of $\Aut(F_n)$ necessarily takes words to words, and since characters of words
are conjugate invariant, $\sigma\in\Out(F_n)$ will take characters to characters.  By Theorem~\ref{thm:Hor},
\[x_i\mapsto \sigma(x_i)\in P_n, \] so that $\sigma$ induces an element of $\Aut_{I_n}(P_n)$ (see
Whittemore~\cite{Wh} or Horowitz~\cite{Hor} for a discussion on this).  Thus $\sigma$ determines a polynomial map
$\widehat\sigma$ on the coordinates of $\C^{2^n-1}$.  In Section~\ref{sec:OutFn}, it is shown that
$\widehat\sigma$ is an automorphism of $\C^{2^n-1}$.  Furthermore, it is shown there that $\det\left(
Jac(\widehat\sigma)\right)=\pm 1$.  Let $\nu$ be the holomorphic $2^n-1$ form on $\C^{2^n-1}$ given by \[
\nu=\bigwedge_{i=1}^{2^n-1} dx_i\in\bigwedge {}^{2^n-1}\C^{2^n-1}.\] This is the complex volume form formed by the
differentials of the Horowitz generators.  Any automorphism of $\C^{2^n-1}$ will take $\nu$ to a functional
multiple of itself (i.e.,
\[ \nu\longmapsto \det\left( Jac(\widehat\sigma))\right) \cdot \nu.\]  By Proposition~\ref{thm:Jac1}, $\det\left(
Jac(\widehat\sigma))\right) \in \{-1,+1\}$.  Thus, up to the sign of $\nu$, $\nu$ is invariant under
$\widehat\sigma $.
\end{proof}

\begin{theorem}  Let $F_n=\pi_1(S)$ for $S$ a compact surface with boundary.  Then $\MCG(S)$ acts as
volume preserving automorphisms of $V_{F_n,\R}\subset\R^{2^n-1}$ which extend to volume preserving polynomial
automorphisms of $\R^{2^n-1}$.\end{theorem}

If we restrict to the real points, then $V_{F_n,\R}=V_{F_n}\cap\R^{2^n-1}$ is a real algebraic set. The
$\SL(\C^2)$-characters of $F_n$ which are real valued come from one of the real forms of $\SL(\C^2)$.  Morgan and
Shalen \cite{MS} show that the real forms of $\SL(\C^2)$ come in two types: $\SL(\R^2)$ and $\SU(2)$.  Let
$F_n=\pi_1(S)$ for some compact orientable surface $S$.  Then it follows that $S$ is of genus-$g$ with
$k=n-2g+1>0$ boundary components.  Here $V_{F_n}$ has the additional structure of a Poisson variety, whose
symplectic leaves are the inverse images of points of the Casimir map:  Let $\{C_i\}_{i=1}^k \subset F_n=\pi_1(S)$
be the set of simple loops homotopic to each of the boundary components of $S$.  Then the map
\[ F:\Hom(S,\SL(\C))\longrightarrow \C^k\] where $F(\phi)= (tr_{C_1}(\phi),\ldots,tr_{C_k}(\phi))=(c_1,\ldots,c_k)$ is conjugation invariant
and descends to a map on $V_{F_n}$ (Recall that $c_i=tr_{C_i}$).  The Poisson structure defines a symplectic
structure on these leaves.  It is known that the mapping class group of the surface $MCG(S)\subset\Out(pi_1(S))$
preserves this Poisson structure, and acts symplectically on each of the leaves.  Hence, the symplectic volume
$\nu_s$ on each of these leaves is preserved under the action of mapping classes.  One can define a volume form
$\nu_\ell$ from the leaf space on which the mapping class group acts trivially (it preserves the leaf structure),
Then, it is easy to see that for $\sigma\in\MCG(S)$, $\sigma$ preserves the volume form on $V_{F_n}$ given by
$\nu_s\wedge F^*\nu_\ell$.

To prove the theorem, for $\sigma\in\Out(\pi_1(S))$, let $\sigma_H\in\Aut(H_1(S;\Z))$ be its associated linear
action on first homology.

\setcounter{section}{6} \setcounter{theorem}{0}

\begin{lemma} For $\sigma\in\MCG(S)$, $\sigma_H$ is unimodular.\end{lemma}

For any particular surface $S$ with $F_n$ as its fundamental group, a presentation of $\MCG(S)$ would lead to a
direct calculation of the corresponding action on the abelianization of $F_n$, which is $H_1(S)$.  Instead, to
prove this assertion without regard to any particular choice of $S$, we will appeal via Poincar\'e duality to the
effect of a mapping class on first cohomology.  We thank Bill Goldman for this suggestion, and a discussion of the
proof.

\begin{proof} To start, let $\partial S=\emptyset$ (while $\pi_1(S)$ would not be free in this case, this case
is instructive).  $H_1(S)$ is Poincar\'e dual to $H^1(S)$.  There is a nondegenerate, skew-symmetric, bilinear
2-form on $H^1(S)$ given by cup product.  \[ \omega:H^1(S)\times H^1(S)\longrightarrow H^2(S)\cong\Z.\]  This
structure forms the basis for the construction of a symplectic structure on the $G$-character variety of $S$, for
$G$ a Lie group with an orthogonal structure (see Goldman~\cite{Go}).  It is known that the induced map $\sigma^*$
preserves this cup product, and acts symplectically on the $G$-character variety of $S$.  Hence the linear action
on $H^1(S)$ and thus on $H_1(S)$ is by symplectic matrix, which is unimodular.

For our case, where necessarily $\partial S\not=\emptyset$, $H^2(S)=0$, the above construction fails.  There is a
corresponding relative version of the cup product, given by \[ \omega: H^1(S,\partial S)\times H^1(S,\partial
S)\longrightarrow H^2(S,\partial S)\cong\Z,\] which is degenerate.  However, by restricting this product to the
parabolic cohomology \[ \omega_p: H_p^1(S,\partial S)\times H_p^1(S,\partial S)\longrightarrow H_p^2(S,\partial
S)\cong H^2(S,\partial S),\] the 2-form is nondegenerate.  Note that parabolic cocycles relative to $\partial S$
are relative cocycles of $S$ which restrict to coboundaries on $\partial S$ (see Weil~\cite{We}).  This is
essentially the construction of Huebschmann~\cite{Hueb} to establish the symplectic structure on the symplectic
leaves of the Poisson character variety of a surface with boundary mentioned above.  $\sigma^*$ preserves this cup
product, and leaves invariant the short exact sequence \[ 0\longrightarrow H_p^1(S,\partial S)\longrightarrow
H^1(S)\longrightarrow H^1(\partial S)\longrightarrow 0.\]  By definition, mapping classes fix pointwise the
components of $\partial S$.  Hence, the determinant of $\sigma^*$ on $H^1(S)$ must equal the product of the
determinant of $\sigma^*$ on $H_p^1(S,\partial S)$ and that of $\sigma^*$ on $H^1(\partial S)$.  As mapping
classes of the surface act identically on the boundary cohomology, $\sigma^*$ on $H^1(S)$ must be
unimodular.\end{proof}

\begin{lemma} For $\sigma\in\Out(F_n)$, \[ \det(\sigma_H)=\det(Jac(\widehat\sigma)).\]\end{lemma}

The proof of this Lemma can be established via a direct calculation, given the four Nielsen generators of
$\Out(F_n)$ and their associated actions of the relevant spaces.  We leave this for the reader.

\begin{proof}[Proof of Theorem~\ref{thm:SubMain}] The Theorem will now follow as a corollary of the above lemmas.
Since outer automorphisms of $F_n=\pi_1(S)$ which correspond to mapping classes of $S$ lead to unimodular
automorphisms of $H_1(S)$, they will act as polynomial automorphisms of $\R^{2^n-1}$ whose Jacobian is everywhere
1.  By the same argument given in Theorem~\ref{thm:Main} above, an automorphism of $\R^{2^n-1}$ will take the
standard real volume form $\nu$ to a functional multiple of itself, so that \[ \nu\longmapsto \det\left(
Jac(\widehat\sigma)\right) \cdot \nu.\] Since for $\sigma\in\MCG(S)$, we have $\det\left(
Jac(\widehat\sigma)\right) = 1$ everywhere, it follows that $\widehat\sigma$ acts as a volume preserving
automorphism of $\R^{2^n-1}$.
\end{proof}

\section{Examples}\label{sec:Examples}

\begin{example} Let $F_2=\langle A,B\rangle$, ordered so that $A>B$.  The Horowitz generating set is given by the
three characters $a$, $b$, and $ab$.  It was proven by Fricke~\cite{Fr} that the character ring of $F_2$ is freely
generated by these three functions, and $I_2$ is trivial.  Hence $V_{F_2}\equiv\C^3$ (compare also
Horowitz~\cite{Hor1} and Whittemore~\cite{Wh}).  $\Out(F_2)$ acts as polynomial automorphisms of $\C^3$ via the
four generators

\small
\[
\widehat\Phi_1:\begin{array}{rcl}
a&\mapsto&ab\\
b&\mapsto&b\\
ab&\mapsto&b\cdot ab-a\end{array} \quad \widehat\Phi_2,\widehat\Phi_3:\begin{array}{rcl}
a&\mapsto&b\\
b&\mapsto&a\\
ab&\mapsto&ab\end{array} \quad \widehat\Phi_4:\begin{array}{rcl}
a&\mapsto&a\\
b&\mapsto&b\\
ab&\mapsto&a\cdot b-ab\end{array}\]

\normalsize
\end{example}

\begin{example} For $F_3=\langle A,B,C\rangle$, where $A>B>C$, the Fricke character ring is the quotient of
\[P_3=\Z[a,b,c,ab,ac,bc,abc]\] cut out by the single generator of $I_3$ given by the Fricke Relation: \small
\begin{eqnarray*} p_{abc}=p_{y_1}&=&a\cdot b\cdot c\cdot abc-a\cdot b\cdot ab-a\cdot c\cdot ac+a\cdot bc\cdot
abc-b\cdot c\cdot bc+b\cdot ac\cdot abc\\ &&+c\cdot ab\cdot abc+ab\cdot ac\cdot
bc+a^2+b^2+c^2+ab^2+ac^2+bc^2+abc^2-4.\end{eqnarray*}\normalsize Hence $V_{F_n}\subset\C^7$ is the zero locus of
$p_{y_1}$.  In this case, The action of $\Out(F_n)$ is via the four generators \tiny
\[
\begin{array}{rcl}&\widehat\Phi_1 &\\ \\
a&\mapsto&ab\\
b&\mapsto&b\\
c&\mapsto&c\\
ab&\mapsto&b\cdot ab-a\\
ac&\mapsto&abc\\
bc&\mapsto&bc\\
abc&\mapsto&b\cdot abc-ac\end{array} \hskip -.05truein\quad
\begin{array}{rcl}&\widehat\Phi_2 &\\ \\
a&\mapsto&b\\
b&\mapsto&a\\
c&\mapsto&c\\
ab&\mapsto&ab\\
ac&\mapsto&bc\\
bc&\mapsto&ac\\
abc&\mapsto&-a\cdot b\cdot c+a\cdot bc\\ &&+b\cdot ac+c\cdot ab-abc\end{array}\hskip -.05truein\quad
\begin{array}{rcl}&\widehat\Phi_3 &\\ \\
a&\mapsto&b\\
b&\mapsto&c\\
c&\mapsto&a\\
ab&\mapsto&bc\\
ac&\mapsto&ab\\
bc&\mapsto&ac\\
abc&\mapsto&abc\end{array}\hskip -.05truein\quad
\begin{array}{rcl}&\widehat\Phi_4 &\\ \\
a&\mapsto&a\\
b&\mapsto&b\\
c&\mapsto&c\\
ab&\mapsto&a\cdot b-ab\\
ac&\mapsto&a\cdot c-ac\\
bc&\mapsto&bc\\
abc&\mapsto&a\cdot bc-abc\end{array}\]

\normalsize

\end{example}

\begin{example}\label{ex:F4} Let $F_4=\langle A,B,C,D\rangle$, ordered so that $A>B>C>D$.  Then the Horowitz generating set is
\begin{eqnarray*} \{x_i\}_{i=1}^{15}&=&\{a,b,c,d,ab,ac,ad,bc,bd,cd,abc,abd,acd,bcd,abcd\} \\
&=&\{l,m,n,o,p,q,r,s,t,u,v,w,x,y,z\},\end{eqnarray*} where in the last equality, we use a notation a bit more
amenable to the number of variables.  Notice here that \[ \{x_j\}_{j=10}^{15} = \{y_i\}_{i=1}^6 =
\{u,v,w,x,y,z\}.\]  $V_{F_4}\subset\C^{15}$ is cut out via the 6 polynomials (note the first is given according to
Case 1 above, and the others are Case 2):

\tiny
\begin{eqnarray*} p_{y_1}&=& l^2m^2n^2o^2-l^2m^2nou-l^2mn^2ot-l^2mno^2s-lm^2n^2or-lm^2no^2q
-lmn^2o^2p-l^2m^2n^2\\ &&-l^2m^2o^2+l^2mntu+l^2mosu-l^2n^2o^2+l^2nost +lm^2nru+lm^2oqu+lmn^2rt+lmo^2qs\\
&&+ln^2opt+lno^2ps-m^2n^2o^2+m^2noqr+mn^2opr+mno^2pq+2l^2mns+2l^2mot+2l^2nou
\\ &&-l^2stu+2lm^2nq+2lm^2or+2lmn^2p+2lmo^2p+lmpu^2
-lmqtu-lmrsu+2ln^2or+2lno^2q\\ &&-lnptu+lnqt^2-lnrst -lopsu-loqst+lors^2+2m^2nou-m^2qru+2mn^2ot+2mno^2s\\ &&
-mnpru-mnqrt+mnr^2s-mopqu+moq^2t-moqrs-n^2prt+nop^2u -nopqt-noprs\\ &&-o^2pqs-l^2s^2-l^2t^2-l^2u^2-2lmqs-2lmrt
-2lnps-2lnru-2lopt-2loqu-m^2q^2-m^2r^2\\ &&-m^2u^2-2mnpq
-2mntu-2mopr-2mosu-n^2p^2-n^2r^2-n^2t^2-2noqr-2nost-o^2p^2
\\ &&-o^2q^2-o^2s^2-p^2u^2+2pqtu+2prsu-q^2t^2+2qrst-r^2s^2-4lmp
-4lnq-4lor-4mns\\ &&-4mot-4nou+4pqs+4prt+4qru+4stu+4l^2+4m^2+4n^2+4o^2+4p^2+4q^2+4r^2\\ &&+4s^2+4t^2+4u^2-16
\\ p_{y_2}&=&lmnv-lmp-lnq-lsv-mns-mqv-npv+pqs
+l^2+m^2+n^2+p^2+q^2+s^2+v^2-4
\\ p_{y_3}&=&lmow-lmp-lor-ltw-mot-mrw-opw+prt
+l^2+m^2+o^2+p^2+r^2+t^2+w^2-4
\\ p_{y_4}&=&lnox-lnq-lor-lux-nou-nrx-oqx+qru
+l^2+n^2+o^2+q^2+r^2+u^2+x^2-4
\\ p_{y_5}&=&mnoy-mns-mot-muy-nou-nty-osy+stu
+m^2+n^2+o^2+s^2+t^2+u^2+y^2-4
\\ p_{y_6}&=&lmuz-lmp-lux-lyz-muy-mxz-puz+pxy+l^2+m^2+p^2+u^2+x^2+y^2+z^2-4\end{eqnarray*} \normalsize

As automorphisms of $\C^{15}$ then, we have \tiny
\[
\widehat\Phi_1:\begin{array}{rcl}
l&\mapsto&p\\
m&\mapsto&m\\
n&\mapsto&n\\
o&\mapsto&o\\
p&\mapsto&mp-l\\
q&\mapsto&v\\
r&\mapsto&w\\
s&\mapsto&s\\
t&\mapsto&t\\
u&\mapsto&u\\
v&\mapsto&mv-q\\
w&\mapsto&mw-r\\
x&\mapsto&z\\
y&\mapsto&y\\
z&\mapsto&mz-x\end{array}\hskip -.1truein \widehat\Phi_2:\begin{array}{rcl}
l&\mapsto&m\\
m&\mapsto&l\\
n&\mapsto&n\\
o&\mapsto&o\\
p&\mapsto&p\\
q&\mapsto&s\\
r&\mapsto&t\\
s&\mapsto&q\\
t&\mapsto&r\\
u&\mapsto&u\\
v&\mapsto&-lmn+ls+mq+np-v\\
w&\mapsto&-lmo+lt+mr+op-w\\
x&\mapsto&y\\
y&\mapsto&x\\
z&\mapsto&-lmu+ly+mx+pu-z\end{array}\hskip -.15truein \widehat\Phi_3:\begin{array}{rcl}
l&\mapsto&m\\
m&\mapsto&n\\
n&\mapsto&o\\
o&\mapsto&l\\
p&\mapsto&s\\
q&\mapsto&t\\
r&\mapsto&p\\
s&\mapsto&u\\
t&\mapsto&q\\
u&\mapsto&r\\
v&\mapsto&y\\
w&\mapsto&v\\
x&\mapsto&w\\
y&\mapsto&x\\
z&\mapsto&z\end{array}\widehat\Phi_4:\begin{array}{rcl}
l&\mapsto&l\\
m&\mapsto&m\\
n&\mapsto&n\\
o&\mapsto&o\\
p&\mapsto&lm-p\\
q&\mapsto&ln-q\\
r&\mapsto&lo-r\\
s&\mapsto&s\\
t&\mapsto&t\\
u&\mapsto&u\\
v&\mapsto&ls-v\\
w&\mapsto&lt-w\\
x&\mapsto&lu-x\\
y&\mapsto&y\\
z&\mapsto&ly-z\end{array}\]

\normalsize

One can easily verify by direct computation that $\widehat\Phi_2$ and $\widehat\Phi_4$ are involutions,
$\widehat\Phi_3$ is of order $n=4$, and that the Jacobian of any of them has constant determinant $\pm 1$.
\end{example}

\begin{example}\label{ex:GAMA} In~\cite{GAMA}, Gonz\'alez-Acu\~na and Montesinos-Amilibia build a set of
polynomials to generate $I_n$ using only the Horowitz generators whose basic words are of length three or less.
These $p=\frac{n(n^2+5)}{6}$ indeterminates provide a complete set of generators, and any Horowitz generator of
word length four or more can be written as a polynomial of the others.  Indeed, for $n=4$, using the above
variables in Example~\ref{ex:F4} above, the Horowitz generator $abcd$ is identified with $z$.  $z$ can be
eliminated through the identity (\cite{GAMA}, Lemma 4.1.1) \[2z=lmno-lmu-los-mnr-nop+ly+mx+nw+ov-qt+rs.\]  The
induced automorphism $\widehat\Phi_1$ is the polynomial map\small \[ \widehat\Phi_1:\begin{array}{rcl}
l&\mapsto&p\\
m&\mapsto&m\\
n&\mapsto&n\\
o&\mapsto&o\\
p&\mapsto&mp-l\\
q&\mapsto&v\\
r&\mapsto&w\\
s&\mapsto&s\\
t&\mapsto&t\\
u&\mapsto&u\\
v&\mapsto&mv-q\\
w&\mapsto&mw-r\\
x&\mapsto&z = \frac{1}{2}(lmno-lmu-los-mnr-nop+ly+mx+nw+ov-qt+rs)\\
y&\mapsto&y\end{array}.\]\normalsize

A quick calculation yields that \[ \det\left( Jac(\widehat\Phi_1)\right) = \frac{1}{2} m,\] which vanishes
wherever $m$ does.\end{example}



\makeatletter \renewcommand{\@biblabel}[1]{\hfill#1.}\makeatother
\renewcommand{\bysame}{\leavevmode\hbox to3em{\hrulefill}\,}


\begin{thebibliography}{20}

\bibitem{Fr}    Fricke, R., and Klein, F.,
        {\em Vorlesungen \"uber die Theorie der automorphem Functionen,\/} Vol. 1, pp. 365-370.  Leipzig: B.G. Teubner              
1897.  Reprint:  New York Juhnson Reprint Corporation (Academic Press) 1965.

\bibitem{GAMA}  Gonz\'alez-Acu\~na, F., and Montesinos-Amilibia, J.M.,
        {\em On the character variety of group representations in $\SL(2,C)$ and $\PSL(2,C)$. ,\/} Math. Z. {\bf 214}
(1993), no. 4, 627-652.

\bibitem{Go}    Goldman, W.,
        {\em The symplectic nature of the fundamental groups of surfaces,\/} Adv. Math. {\bf 54} (1984), 200-225.

\bibitem{Hor1}  Horowitz, R.,
        {\em Characters of free groups represented in the two-dimensional special linear group,\/}
        Comm. Pure Appl. Math. {\bf 25} (1972), 635-649.

\bibitem{Hor}   Horowitz, R.,
        {\em Induced automorphisms on Fricke characters of free groups,\/} Trans. AMS {\bf 208} (1975), 41-50.

\bibitem{Hueb}  Huebschmann, J.,
        {\em On the variation of the Poisson structures of certain moduli spaces,\/} Math. Ann. {\bf 319} (2001), 267-310.

\bibitem{Mag}   Magnus, W.,
        {\em Rings of Fricke characters and automorphism groups of free groups,\/} Math. Z. {\bf 170} (1980), 91-102.

\bibitem{Mc}    McCool, J.,
        {\em Two-dimensional linear characters and automorphisms of free groups,\/}
        Combinatorial group theory (College Park, MD, 1988), 131-138,
        Contemp. Math. {\bf 109}, Amer. Math. Soc., Providence, RI, 1990.

\bibitem{MS}    Morgan, J., and Shalen, P.,
        {\em Valuations, trees, and degenerations of hyperbolic structures I ,\/} Ann. Math. {\bf 120}
        (1984), 401-476.

\bibitem{Niel}  Nielsen, J.,
        {\em Die Isomorphismengruppe der freien Gruppen\/} Math. Ann. {\bf 91} (1924), 169-209.

\bibitem{Th}    Thurston, W.,
        {\em On the geometry and dynamics of diffeomorphisms of surfaces,\/} Bull. AMS {\bf 19} (1988), 417-431.

\bibitem{Wang}  Wang, S.S.-S.,
        {\em A Jacobian criterion for separability,\/} J. Algebra {\bf 65} (1980), 453-440.

\bibitem{We}    Weil, A.,
        {\em Remarks on the cohomology of groups,\/} Ann. Math.(2) {\bf 80} (1964), 149-157.

\bibitem{Wh}    Whittemore, A.,
        {\em On special linear characters of free groups of rank $n\le 4$,\/} Proc. AMS {\bf 40} (1973), 383-388.

\end{thebibliography}
\end{document}